# Consistency of the $\alpha$-trimming of a probability. Applications to central regions

IGNACIO CASCOS[1] and MIGUEL LÓPEZ-DÍAZ[2]

[1]*Departamento de Estadística, Universidad Carlos III de Madrid, Av. Universidad 30, E-28911 Leganés (Madrid), Spain. E-mail: ignacio.cascos@uc3m.es*
[2]*Departamento de Estadística e I.O. y D.M., Universidad de Oviedo Facultad de Ciencias, c/ Calvo Sotelo s/n, E-33007 Oviedo, Spain. E-mail: mld@uniovi.es*

The sequence of $\alpha$-trimmings of empirical probabilities is shown to converge, in the Painlevé–Kuratowski sense, on the class of probability measures endowed with the weak topology, to the $\alpha$-trimming of the population probability. Such a result is applied to the study of the asymptotic behaviour of central regions based on the trimming of a probability.

## 1. Introduction

Depth-trimmed regions are sets of central points with respect to a multivariate probability distribution or a data cloud. Given a mapping that associates each point with its degree of centrality, which is traditionally called *depth*, the set of all points whose depth is at least a given value constitutes a central or depth-trimmed region. Therefore, central regions are set-valued location estimates. Their applications range from simple descriptive statistics to multivariate generalizations of univariate quantile-based techniques. Numerous proposals of depth-trimmed regions have been considered in the statistical literature in recent years; see, for instance, Massé and Theodorescu [8] for half-space trimmed regions, Koshevoy and Mosler [7] for zonoid trimmed regions, Cascos and López-Díaz [3] for integral trimmed regions and Zuo and Serfling [12] for a unifying treatment of several notions of central regions.

Cascos and López-Díaz [3] introduce the integral trimmed regions as a family of depth-trimmed regions that is generated by a set (or family of sets) of functions. Several classical families of depth-trimmed regions like the zonoid or the halfspace regions can be obtained in the framework of integral trimmed regions for suitable families of functions. Integral trimmed regions have been used, among other families of depth-trimmed regions, by Cascos and Molchanov [4] to assess the financial risk of a vector portfolio. The concept of the $\alpha$-trimming of a probability plays a crucial role in the construction of integral trimmed regions and thus a comprehensive study of its properties would constitute a significant advance in the theoretical study of depth-trimmed regions. In particular, the analysis of the behaviour of the $\alpha$-trimming of empirical probabilities is necessary for the study of empirical integral trimmed regions.

Motivated by this, the problem of the consistency of the $\alpha$-trimming of a sequence of empirical probabilities is analyzed in this paper. On the basis of such a consistency, applications to integral trimmed regions and to other trimmed regions based on the trimming of a probability will be developed.

Throughout the paper, the symbol $\mathcal{P}$ will stand for the class of probability measures on $(R^d, \mathcal{B}_d)$, where $\mathcal{B}_d$ denotes the Borel $\sigma$-algebra on $R^d$. This class will be endowed with the weak topology.

As has already been mentioned, the main results of the paper are focused on the concept of $\alpha$-trimming of a probability of $\mathcal{P}$, which can be defined as follows.

**Definition 1.** *Given $P \in \mathcal{P}$ and $\alpha \in (0,1]$, the $\alpha$-trimming of $P$, denoted by $P^\alpha$, is given by $P^\alpha = \{Q \in \mathcal{P} : Q \leq \alpha^{-1} P\}$, that is, the set of all probability measures $Q$ on $(R^d, \mathcal{B}_d)$ such that $Q(B) \leq \alpha^{-1} P(B)$ for all $B \in \mathcal{B}_d$.*

This set of probabilities was formally introduced by Cascos and López-Díaz [3] analyzing several properties of such a set. Sets of probabilities that are bounded by a transformation of another probability have been studied by, among others, Huber and Strassen [6] in the field of robust statistics and Artzner *et al.* [1] in the field of financial mathematics.

From now on, given $P \in \mathcal{P}$, $\{X_i\}_i$ will be a sequence of independent and identically distributed random vectors, where $X_1$ induces the probability $P$. Let $\{P_n\}_n \subset \mathcal{P}$ be the sequence of empirical probabilities generated by $\{X_i\}_i$.

Given $P \in \mathcal{P}, \{P_n\}_n$ a sequence of empirical probabilities and $\alpha \in (0,1]$, this paper analyzes the behaviour of the sequence of probability sets $\{P_n^\alpha\}_n$ in relation to the set $P^\alpha$. Namely, it is proved that the sequence $\{P_n^\alpha\}_n$ converges (almost surely) to $P^\alpha$.

The usual concepts of lower limit and upper limit of sequences of subsets of a space with a convergence criterion will be considered. Limits of sequences of sets are understood in the Painlevé–Kuratowski sense, that is, a sequence of sets converges if its lower and upper limit coincide, this being equal to the limit.

Since $\mathcal{P}$ is endowed with the convergence of the weak topology, we have

$$\limsup_n P_n^\alpha = \left\{ Q \in \mathcal{P} : \exists \{Q_{n_k}\}_k, Q_{n_k} \in P_{n_k}^\alpha \ \forall k \in N \text{ with } \lim_k Q_{n_k} = Q \right\}$$

and

$$\liminf_n P_n^\alpha = \left\{ Q \in \mathcal{P} : \exists \{Q_n\}_n \text{ with } Q_n \in P_n^\alpha \ \forall n \in N, \text{ and } \lim_n Q_n = Q \right\},$$

where the limits inside the brackets are considered in the convergence of the weak topology.

The structure of the paper is as follows. Section 2 is devoted to study of the limit behaviour of the $\alpha$-trimming of an empirical probability and applications of this limit behaviour to depth-trimmed regions are derived in Section 3.

## 2. Trimming of empirical probabilities

The aim of this section is to show that given a probability $P \in \mathcal{P}$ and $\alpha \in (0,1]$, the set of probabilities $P^\alpha$ is the limit (almost surely) of the sequence of the sets of probabilities $\{P_n^\alpha\}_n$. More specifically, the main result will conclude that

$$\liminf_n P_n^\alpha = P^\alpha = \limsup_n P_n^\alpha \qquad \text{a.s.}$$

This section contains three parts. First, we analyze some supporting properties in relation to the topological structure of the trimming of a probability. The second part contains some key technical results for the main purpose of the section. Finally, the study of sequences of the trimming of empirical probabilities is developed.

We will now state some results about topological properties of the trimming of a probability.

A family of probabilities $\Pi \subset \mathcal{P}$ is said to be *tight* if, for each $\varepsilon > 0$, there exists a compact set $K \subset R^d$ such that $P(K) > 1 - \varepsilon$ for all $P \in \Pi$. It is well known that if $P \in \mathcal{P}$, then the family $\{P\}$ is tight.

**Lemma 2.** *Let $\Pi \subset \mathcal{P}$ be a tight family and $\alpha \in (0,1]$. Then, $\{Q \in H^\alpha : H \in \Pi\}$ is tight.*

**Proof.** Since $\Pi$ is tight, given $\varepsilon > 0$, there exists a compact set $K \subset R^d$ such that $H(R^d \setminus K) < \varepsilon \alpha$ for all $H \in \Pi$. So, for any $Q \in H^\alpha$, we have that $Q(R^d \setminus K) \leq \alpha^{-1} H(R^d \setminus K) < \varepsilon$, which implies the result. □

**Lemma 3.** *For all $\alpha \in (0,1]$ and $P \in \mathcal{P}$, the family $P^\alpha$ is compact for the weak topology.*

**Proof.** By Lemma 2, the family $P^\alpha$ is tight, so its closure in the weak topology $\mathrm{cl}P^\alpha$ is compact (see, e.g., Billingsley [2]). Thus, it is sufficient to prove that $\mathrm{cl}P^\alpha = P^\alpha$.

Let $Q \in \mathrm{cl}P^\alpha$. There then exists $\{Q_n\}_n \subset P^\alpha$ with $\lim_n Q_n = Q$ in the convergence of the weak topology.

Let $G$ be an open set in $R^d$, so $Q(G) \leq \liminf_n Q_n(G) \leq \alpha^{-1} P(G)$. If $A \in \mathcal{B}_d$, we can take a decreasing sequence $\{G_n\}_n$ of open sets in $R^d$ with $A \subset G_n$ and $P(G_n \setminus A) < 1/n$ (observe that $P$ is regular), so $Q(A) \leq Q(G_n) \leq \alpha^{-1} P(G_n)$ for all $n \in N$. Then, $Q(A) \leq \alpha^{-1} P(A)$, which implies that $Q \in P^\alpha$, so $P^\alpha$ is compact. □

**Lemma 4.** *Let $P \in \mathcal{P}$ and $\alpha \in (0,1]$. There then exists a countable set $\mathcal{D} \subset P^\alpha$ such that $\mathcal{D}$ is dense in $P^\alpha$ in the weak topology.*

**Proof.** Note that the weak topology is metrizable, for instance, by means of the Prohorov metric and $P^\alpha$ is compact in such a topology, which trivially implies the result. □

In order to obtain the main result of this section, some essential technical results are stated below.

The following proposition shows how to construct an element of the class $P_{n+1}^\alpha$ by means of a probability of the class $P_n^\alpha$. This result will play a crucial role in posterior constructions of sequences of probabilities.

**Proposition 5.** *Let $P \in \mathcal{P}$ and $\alpha \in (0,1]$. If $Q_n \in P_n^\alpha$, then the probability $Q_{n+1}$ given by*

$$Q_{n+1} = \frac{n}{n+1} Q_n + \frac{1}{n+1}((n+1)P_{n+1} - nP_n)$$

*belongs to the class $P_{n+1}^\alpha$.*

**Proof.** It should be pointed out that the mapping $(n+1)P_{n+1} - nP_n : \mathcal{B}_d \to R$ satisfies $((n+1)P_{n+1} - nP_n)(B) = I_B(X_{n+1})$ for all $B \in \mathcal{B}_d$ and is therefore a probability. As a consequence, $Q_{n+1}$ belongs to $\mathcal{P}$.

Further, $(n+1)P_{n+1} - nP_n$ is a discrete probability. Since $Q_n \in P_n^\alpha$, it follows that $Q_n$ is absolutely continuous with respect to $P_n$, which implies that $Q_n$ is also discrete and so the same holds for $Q_{n+1}$.

Thus, in order to check that $Q_{n+1} \in P_{n+1}^\alpha$, it is sufficient to prove that $Q_{n+1}(\{x\}) \leq \alpha^{-1} P_{n+1}(\{x\})$ for all $x \in R^d$.

Let $x \in R^d$. If $x \neq X_{n+1}$, then $((n+1)P_{n+1} - nP_n)(\{x\}) = 0$, so

$$Q_{n+1}(\{x\}) = \frac{n}{n+1} Q_n(\{x\}) \leq \alpha^{-1} \frac{n}{n+1} P_n(\{x\}) = \alpha^{-1} P_{n+1}(\{x\}).$$

If $x = X_{n+1}$, then $((n+1)P_{n+1} - nP_n)(\{x\}) = 1$, which implies that

$$Q_{n+1}(\{x\}) = \frac{n}{n+1} Q_n(\{x\}) + \frac{1}{n+1} \leq \alpha^{-1} \frac{n}{n+1} P_n(\{x\}) + \frac{1}{n+1}$$

$$\leq \alpha^{-1} \frac{nP_n(\{x\}) + 1}{n+1} = \alpha^{-1} P_{n+1}(\{x\})$$

and this finishes the proof. □

Note that, given $Q_n \in P_n^\alpha$, the reiteration of the construction of probabilities described in the statement of the above proposition leads to the probability

$$Q_{n+k}(\cdot) = \frac{n}{n+k} Q_n(\cdot) + \frac{1}{n+k}(I_{(\cdot)}(X_{n+1}) + I_{(\cdot)}(X_{n+2}) + \cdots + I_{(\cdot)}(X_{n+k})) \quad (1)$$

for all $k \in \{1, 2, \ldots\}$.

A technical lemma on stopping times that will be applied in posterior results is stated below.

**Lemma 6.** *Let $\{Z_i\}_i$ be a sequence of independent and identically distributed random variables with $EZ_1 = 0$. For all $n \in N$, let $S_n = \sum_1^n Z_i$. We define $Y_1 = \min\{n \in N : S_n \geq 0\}$ and for $k \geq 1$, let $Y_{k+1} = \min\{n > Y_k : S_n \geq 0\}$.*

*Then,*

$$\lim_k \frac{Y_{k+1} - Y_k}{Y_k} = 0 \qquad a.s.$$

**Proof.** Given $Z_1, \ldots, Z_n$, the smallest $\sigma$-field which makes these mappings measurable will be denoted by $\sigma(Z_1, \ldots, Z_n)$.

In the first place, we should point out that the probability that the event $\sum_1^n Z_i \geq 0$ occurs infinitely often is equal to 1.

We define $W_0 = 0$, $W_1 = \min\{n > 0 : \sum_1^n Z_i \geq 0\}$ and, for $k \geq 1$, let $W_{k+1} = \min\{n > 0 : \sum_{W_k+1}^n Z_i \geq 0\}$.

Note that $W_k$ is a stopping time since, for all $k \in N$, it holds that $W_k$ is $\sigma(Z_1, \ldots, Z_{W_k})$-measurable. Also, $W_k \in N$ a.s. for all $k \in N$.

Moreover, the random variables $\{W_{k+1} - W_k\}_k$ are independent since $W_{k+1} - W_k$ is $\sigma(Z_{W_k+1}, \ldots, Z_{W_{k+1}})$-measurable for all $k \in N$ and is obviously identically distributed.

Well-known results on stopping times state that $EW_1 < \infty$ since $EZ_1 \in R$ (see, e.g., Wijsman [11] or Durrett [5]), so $E(W_{k+1} - W_k) = E(W_1 - W_0) = EW_1 < \infty$.

The strong law of large numbers applied to the sequence $\{W_{k+1} - W_k\}_k$ leads to

$$\lim_n \frac{W_{n+1}}{n} = EW_1 \quad \text{a.s., hence} \quad \lim_n \frac{W_{n+1} - W_n}{n} = 0 \quad \text{a.s.}$$

On the other hand, it can immediately be seen that for each $l \in N$, there exists $k_l \in N$ such that $Y_{l+1} - Y_l \leq W_{k_l+1} - W_{k_l}$.

It is enough to consider $k_l$ with $k_l = \sup\{n : W_n \leq Y_l\}$. Observe that such a $k_l$ exists a.s. since $W_1 = Y_1$. In accordance with the definitions of $\{W_i\}_i, \{Y_i\}_i$ and $k_l$, we obviously have $W_{k_l+1} > Y_l$.

Moreover, $Y_{l+1} - Y_l \leq W_{k_l+1} - W_{k_l}$, otherwise $S_{W_{k_l+1}} \geq 0$, which contradicts the definition of $Y_{l+1}$. In fact, $W_{k_l+1} \geq Y_{l+1}$.

Trivially, $\{Y_i\}_i$ is a non-decreasing sequence, $Y_n \geq n$, $W_n \geq n$ for all $n \in N$ and $\{k_l\}_l$ is nondecreasing, so $\lim_l k_l = \infty$ a.s. Then,

$$0 \leq \limsup_l \frac{Y_{l+1} - Y_l}{Y_l} \leq \limsup_l \frac{W_{k_l+1} - W_{k_l}}{W_{k_l}} \leq \limsup_l \frac{W_{k_l+1} - W_{k_l}}{k_l} = 0 \qquad \text{a.s.,}$$

which proves the lemma. □

Cascos and López-Díaz [3] prove the following result in relation to the trimming of a probability.

**Lemma 7.** *Let $\alpha \in (0, 1]$ and $Q \in P^\alpha$. There then exists a.s. a sequence $\{Q_{n_k}\}_k \subset \mathcal{P}$ such that $Q_{n_k} \in P_{n_k}^\alpha$ for all $k \in N$ and $\lim_k Q_{n_k} = Q$ in the weak topology.*

This result is proved by constructing a particular sequence $\{Q_{n_k}\}_k$ satisfying the above conditions. Since we will make use of such a sequence, it will be described below.

If $Q \in P^\alpha$, then $Q$ is absolutely continuous with respect to $P$. Therefore, there exists a mapping $g: R^d \longrightarrow [0, \infty)$ of class $L^1(P)$ such that for all $B \in \mathcal{B}_d$, it holds that $Q(B) = \int_B g \, dP$, that is, $g$ is a Radon–Nikodym derivative of $Q$ with respect to $P$.

If $g$ is a Radon–Nikodym derivative of $Q$ with respect to $P$, we can explicitly give a sequence $\{Q_{n_k}\}_k$ satisfying the conditions of Lemma 7 as

$$Q_{n_k}(B) = \frac{1}{\|g\|_{L^1(P_{n_k})}} \int_B g \, dP_{n_k} \qquad \text{for all } B \in \mathcal{B}_d, \tag{2}$$

where $n_k$ must satisfy $\|g\|_{L^1(P_{n_k})} \geq 1$.

That is, if we denote by $N_Q$ the set of indexes of the above subsequence, then

$$n \in N_Q \quad \text{if and only if} \quad \frac{1}{n} \sum_{i=1}^n (|g(X_i)| - 1) \geq 0. \tag{3}$$

Note that the probability that the event $n^{-1} \sum_1^n (|g(X_i)| - 1) \geq 0$ occurs infinitely often is equal to 1.

We should indicate that the a.s. existence of a converging subsequence with limit $Q$ in Lemma 7 depends on $g$, a Radon–Nikodym derivative of $Q$ with respect to $P$ and it thus depends on probability $Q$. For a detailed explanation of this construction, we refer to Cascos and López-Díaz [3].

**Proposition 8.** *Let $\{X_i\}_i$ be a sequence of independent and identically distributed random vectors, where $X_1$ induces the probability $P$. Let $g$ be a Radon–Nikodym derivative of $Q$ with respect to $P$. Let $\{n_k\}_k$ be defined by $n_1 = \min\{n \in N : n^{-1} \sum_1^n (|g(X_i)| - 1) \geq 0\}$ and for $k \geq 1$, let $n_{k+1} = \min\{n > n_k : n^{-1} \sum_1^n (|g(X_i)| - 1) \geq 0\}$. It then holds that*

$$\lim_k \frac{n_{k+1} - n_k}{n_k} = 0 \qquad a.s.$$

**Proof.** Note that the random variables $\{|g(X_i)| - 1\}_i$ are independent and identically distributed. Since $E|g(X_i)| = \int_{R^d} g \, dP = Q(R) = 1$, their expected values are equal to 0, thus the result is a direct consequence of Lemma 6. □

It is now proved that any probability of the class $P^\alpha$ belongs almost surely to the inferior limit of the sequence $\{P_n^\alpha\}_n$.

**Proposition 9.** *Let $P \in \mathcal{P}$, $\alpha \in (0, 1]$ and $Q \in P^\alpha$. Then.*

$$Q \in \liminf_n P_n^\alpha \qquad a.s.$$

**Proof.** Lemma 7 states that $Q \in \limsup P_n^\alpha$ a.s. In fact, we have a specific sequence $\{Q_{n_k}\}_k$ with $Q_{n_k} \in P_{n_k}^\alpha$ for all $k \in N$ such that $\lim_k Q_{n_k} = Q$ a.s. in the weak topology.

We recall that $N_Q$ is the set of indices of this sequence which was already characterized in (3).

We will complete this subsequence in the following way. Given $n_1 \in N_Q$, if $n_1 + 1 \notin N_Q$, we define the probability $Q_{n_1+1} \in P^\alpha_{n_1+1}$ by means of Proposition 5, that is,

$$Q_{n_1+1} = \frac{n_1}{n_1+1} Q_{n_1} + \frac{1}{n_1+1}((n_1+1)P_{n_1+1} - n_1 P_{n_1}),$$

continue until $n_1 + k \in N_Q$ and let $n_1 + k = n_2$ (note that $Q_{n_2}$ has already been defined). We then reiterate this process with $n_2$ and so forth with successive indexes.

In this way, for all $n \geq n_1$, we have a probability $Q_n \in P^\alpha_n$ such that if $l \in N_Q$, then $Q_l$ belongs to the subsequence described in (2). In order to finish the proof, we must show that $\lim_n Q_n = Q$ a.s in the weak topology.

Since $\lim_n P_n = P$ a.s. in the weak topology, $\{P_1, P_2, \ldots\}$ is tight a.s. (see, e.g., Billingsley [2]) and so Lemma 2 implies that the family $\bigcup_{n=1}^\infty P^\alpha_n$ satisfies the same property. Hence $\{Q_n\}_{n \geq n_1} \subset \bigcup_{n=n_1}^\infty P^\alpha_n$ is also tight a.s.

As a consequence, in order to prove that $\lim_n Q_n = Q$ a.s., it is enough to show that a.s., given any subsequence $\{Q_{l_k}\}_k \subset \{Q_n\}_n$ with $\lim_k Q_{l_k} = Q'$ for some $Q' \in \mathcal{P}$, we have $Q' = Q$.

Let us consider a subsequence satisfying the above conditions and let $N_{Q'}$ be the set of indexes of such a subsequence.

Clearly, if $|N_{Q'} \cap N_Q| = \infty$, the uniqueness of the limit implies that $Q' = Q$.

Suppose now that $|N_{Q'} \cap N_Q| < \infty$. Then, given $n_k \in N_Q$, we define (whenever it exists) $l'_k = \inf\{l_k \in N_{Q'} : n_k < l_k < n_{k+1}\}$ (note that there exist infinite values of $l'_k$). Observe that $\{Q_{l'_k}\}_k \subset \{Q_{l_k}\}_k$, and so $\lim_k Q_{l'_k} = Q'$.

In accordance with (1), it is trivial that

$$Q_{l'_k}(\cdot) = \frac{n_k}{l'_k} Q_{n_k}(\cdot) + \frac{1}{l'_k}(I_{(\cdot)}(X_{n_k+1}) + I_{(\cdot)}(X_{n_k+2}) + \cdots + I_{(\cdot)}(X_{l'_k})).$$

Moreover,

$$0 \leq \frac{1}{l'_k}(I_{(\cdot)}(X_{n_k+1}) + I_{(\cdot)}(X_{n_k+2}) + \cdots + I_{(\cdot)}(X_{l'_k})) \leq \frac{l'_k - n_k}{l'_k} \leq \frac{n_{k+1} - n_k}{n_k}.$$

Proposition 8 then implies that

$$\lim_k \frac{1}{l'_k}(I_{(\cdot)}(X_{n_k+1}) + I_{(\cdot)}(X_{n_k+2}) + \cdots + I_{(\cdot)}(X_{l'_k})) = 0 \quad \text{a.s.}$$

On the other hand,

$$\lim_k \frac{n_k}{l'_k} = \lim_k \frac{n_k - l'_k}{l'_k} + 1 = 1 \quad \text{a.s.}$$

and so

$$\lim_k Q_{l'_k} = \lim_k Q_{n_k} = Q \quad \text{a.s.}$$

in the weak topology. Then, $Q = Q'$ a.s., which proves the theorem. □

The $\alpha$-trimming of $P$ is (almost surely) contained in the lower limit of the $\alpha$-trimmings of the empirical probabilities, as shown below.

**Theorem 10.** *Let $P \in \mathcal{P}$ and $\alpha \in (0, 1]$. Then,*

$$P^\alpha \subset \liminf_n P_n^\alpha \qquad a.s.$$

**Proof.** Lemma 4 proves that there exists a countable set $\mathcal{D} \subset P^\alpha$ such that $\mathcal{D}$ is dense in $P^\alpha$ in the weak topology. Given any $Q \in \mathcal{D}$, Proposition 9 implies that $Q \in \liminf_n P_n^\alpha$ a.s. and since $\mathcal{D}$ is countable, $\mathcal{D} \subset \liminf_n P_n^\alpha$ a.s.

Finally, $\liminf_n P_n^\alpha$ is closed in the weak topology, thus $\mathrm{cl}\mathcal{D} \subset \liminf_n P_n^\alpha$ a.s., where cl stands for the closure in such a topology.

Since $P^\alpha \subset \mathrm{cl}\mathcal{D}$, we have $P^\alpha \subset \liminf_n P_n^\alpha$ a.s. and this trivially concludes the proof. □

The next result refers to the upper limit of the sequence $\{P_n^\alpha\}_n$.

**Theorem 11.** *Let $P \in \mathcal{P}$ and $\alpha \in (0, 1]$. Then,*

$$\limsup_n P_n^\alpha \subset P^\alpha \qquad a.s.$$

**Proof.** Let $Q \in \limsup_n P_n^\alpha$. There then exists a sequence $\{Q_{n_k}\}_k \subset \mathcal{P}$ with $Q_{n_k} \in P_{n_k}^\alpha$ for all $k \in N$ such that $\lim_k Q_{n_k} = Q$ in the weak topology.

Since $\lim_n P_n = P$ a.s. in the weak topology, for any open set $G$ in $R^d$,

$$Q(G) \leq \liminf_k Q_{n_k}(G) \leq \liminf_k \alpha^{-1} P_{n_k}(G)$$

$$\leq \limsup_k \alpha^{-1} P_{n_k}(\mathrm{cl}G) \leq \alpha^{-1} P(\mathrm{cl}G) \qquad \text{a.s.,}$$

where $\mathrm{cl}G$ stands for the usual closure of the set $G$.

If $F$ is any closed set in $R^d$, consider $G_m = \{x \in R^d : d(x, F) < 1/m\}$, where $d(x, F) = \inf\{\|x - y\| : y \in F\}$ and $m \in N$.

For all $m \in N$, it then holds that $Q(F) \leq Q(G_m) \leq \alpha^{-1} P(\mathrm{cl}G_m)$ a.s. and so $Q(F) \leq \alpha^{-1} P(F)$ a.s.

Since this inequality holds for all closed sets in $R^d$, given $A \in \mathcal{B}_d$ we have only to consider an increasing sequence $\{F_m\}_m$ of closed sets in $R^d$ with $F_m \subset A$, $P(A \setminus F_m) < 1/m$ and $Q(A \setminus F_m) < 1/m$ for all $m \in N$. Note that this is always possible since $P$ and $Q$ are regular probabilities. Therefore, $Q(A) \leq \alpha^{-1} P(A)$ for all $A \in \mathcal{B}_d$ a.s. and so $Q \in P^\alpha$ a.s. □

Finally, the main result of this section is obtained as a consequence of Theorems 10 and 11.

**Theorem 12.** *Let $P \in \mathcal{P}$ and $\{P_n\}_n \subset \mathcal{P}$ be a sequence of empirical probabilities of P. Then, for all $\alpha \in (0,1]$, it holds that the sequence $\{P_n^\alpha\}_n$ converges a.s. in the Painlevé–Kuratowski sense to the set $P^\alpha$, that is,*

$$\lim_n P_n^\alpha = P^\alpha \qquad a.s.$$

## 3. Applications to depth-trimmed regions

Applications of the preceding results to depth-trimmed regions are developed in this section. In the first place, we introduce a new family of depth-trimmed regions based on the trimming of a probability and study its empirical counterpart by means of the results in Section 2. The analysis of empirical integral trimmed regions is then developed.

For each centrality degree $\alpha$, a depth-trimmed region of level $\alpha$ is a set-valued location estimate. Given any (point-valued) location estimate associated with a probability measure, we will build $\alpha$-trimmed regions for each probability $P \in \mathcal{P}$ as the set of all location estimates of the probabilities in the $\alpha$-trimming of $P$. Formally, these regions, which we will call *location trimmed regions*, can be defined as follows.

**Definition 13.** *Let $\alpha \in (0,1]$, let L be a location estimate on the class $\mathcal{P}$ and let $P \in \mathcal{P}$. The location $\alpha$-trimmed region of P induced by L, denoted by $\mathrm{D}_L^\alpha(P)$, is defined as*

$$\mathrm{D}_L^\alpha(P) = \{L(Q) : Q \in P^\alpha\}.$$

Note that some well-established trimmed regions can be obtained by means of location trimmed regions, taking appropriate location estimates.

For instance, consider the zonoid trimmed regions defined by Koshevoy and Mosler [7] as

$$\mathrm{ZD}^\alpha(P) = \left\{ \int x g(x)\,\mathrm{d}P, g \colon R^d \to [0, \alpha^{-1}] \text{ measurable with } \int g(x)\,\mathrm{d}P = 1 \right\}$$

for any $\alpha \in (0,1]$ and $P \in \mathcal{P}$.

It is easy to see that if we consider the expected value as location estimate, that is, $L$ stands for the expected value, then it holds that $\mathrm{ZD}^\alpha(P) = \mathrm{D}_L^\alpha(P)$.

The following proposition concerns, under mild conditions, the consistency of empirical location trimmed regions.

**Theorem 14.** *Let L be a location estimate that is continuous with respect to the weak topology. Then, for any $P \in \mathcal{P}$ and $\alpha \in (0,1]$, the empirical location trimmed regions tend, in the Painlevé–Kuratowski sense, to the population trimmed region, that is,*

$$\mathrm{D}_L^\alpha(P) = \lim_n \mathrm{D}_L^\alpha(P_n) \qquad a.s.$$

**Proof.** Let $x \in D_L^\alpha(P)$. There then exists $Q \in P^\alpha$ such that $x = L(Q)$. Theorem 12 states that $Q \in \lim_n P_n^\alpha$ a.s., so there exists a.s. a sequence $\{Q_n\}_n \subset \mathcal{P}$ with $Q_n \in P_n^\alpha$ for all $n \in N$ such that $\lim_n Q_n = Q$ in the convergence of the weak topology.

Since $L$ is continuous when such a topology is considered, then $x = L(Q) = \lim_n L(Q_n)$ and thus $D_L^\alpha(P) \subset \liminf_n D_L^\alpha(P_n)$ a.s.

Conversely, let $x \in \limsup_n D_L^\alpha(P_n)$. There then exists a sequence $\{x_{n_k}\}_k$ with $x_{n_k} \in D_L^\alpha(P_{n_k})$ such that $x = \lim_k x_{n_k}$. Therefore, there exists $Q_{n_k} \in P_{n_k}^\alpha$ with $x_{n_k} = L(Q_{n_k})$ for all $k \in N$.

Since $\{Q_{n_k}\}_k \subset \bigcup_{n=1}^\infty P_n^\alpha$ and this set is tight, there exists a subsequence $\{Q_{n_{k_l}}\}_l \subset \{Q_{n_k}\}_k$ converging in the weak topology to a certain $Q \in \mathcal{P}$. Theorem 11 implies that $Q \in P^\alpha$ a.s., hence $x = \lim_l x_{n_{k_l}} = \lim_l L(Q_{n_{k_l}}) = L(Q) \in D_L^\alpha(P)$ a.s., which proves that $\limsup_n D_L^\alpha(P_n) \subset D_L^\alpha(P)$ a.s., so the result holds. □

The definition of the integral trimmed regions is based on the concept of $\alpha$-trimming of a probability. Formally, they are defined as follows.

**Definition 15.** *Given $\alpha \in (0,1]$ and a set of measurable functions $\mathcal{F}$ from $R^d$ into $R$, the integral $\alpha$-trimmed region of $P$ with respect to the set $\mathcal{F}$, denoted by $D_\mathcal{F}^\alpha(P)$, is defined as*

$$D_\mathcal{F}^\alpha(P) = \left\{ x \in R^d : \exists Q_x \in P^\alpha \text{ with } f(x) \leq \int f \, dQ_x \text{ for all } f \in \mathcal{F} \right\}$$
$$= \bigcup_{Q \in P^\alpha} \bigcap_{f \in \mathcal{F}} f^{-1}\left(\left(-\infty, \int f \, dQ\right]\right).$$

Integral trimmed regions satisfy, under mild conditions in the generating families of functions, some convenient properties for general trimmed regions, such as being nested, affine invariant, closed or bounded (see Cascos and López-Díaz [3]).

We will obtain a result for the empirical version of integral trimmed regions.

Some concepts which appear in this study are briefly described, in particular the concepts of a Glivenko–Cantelli class for a probability and the concept of an envelope of a family of functions.

Given $P \in \mathcal{P}$, a family of measurable functions (from $R^d$ into $R$) $\mathcal{F} \subset L^1(P)$ is said to be a *Glivenko–Cantelli class* for $P$ if

$$\lim_n \left( \sup_{f \in \mathcal{F}} \left| \int f \, dP_n - \int f \, dP \right| \right) = 0 \quad \text{a.s.}$$

A measurable mapping $f_\mathcal{F}$ is said to be an *envelope* of the family of functions $\mathcal{F}$ if $|f| \leq f_\mathcal{F}$ for all $f \in \mathcal{F}$.

For a more detailed explanation on these concepts and their applications see, for instance, Talagrand [9] and van der Vaart [10].

**Theorem 16.** *If $\mathcal{F}$ is a Glivenko–Cantelli class for $P$ of continuous and bounded functions with an envelope $f_{\mathcal{F}} \in L^1(P)$, then for all $\alpha \in (0,1]$, it holds that*

$$\bigcap_{\varepsilon > 0} \liminf_{n} \mathrm{D}_{\mathcal{F}}^{\alpha,\varepsilon}(P_n) = \mathrm{D}_{\mathcal{F}}^{\alpha}(P) \qquad a.s.,$$

*where*

$$\mathrm{D}_{\mathcal{F}}^{\alpha,\varepsilon}(P) = \left\{ x \in R^d : \exists Q_x \in P^{\alpha} \text{ with } f(x) - \varepsilon \leq \int f \,\mathrm{d}Q_x \text{ for all } f \in \mathcal{F} \right\}.$$

**Proof.** Theorem 27 in Cascos and López-Díaz [3] states that under the conditions of the theorem it holds that

$$\bigcap_{\varepsilon > 0} \limsup_{n} \mathrm{D}_{\mathcal{F}}^{\alpha,\varepsilon}(P_n) = \mathrm{D}_{\mathcal{F}}^{\alpha}(P) \qquad \text{a.s.}$$

Therefore,

$$\bigcap_{\varepsilon > 0} \liminf_{n} \mathrm{D}_{\mathcal{F}}^{\alpha,\varepsilon}(P_n) \subset \mathrm{D}_{\mathcal{F}}^{\alpha}(P) \qquad \text{a.s.}$$

and only the other inclusion remains to be proven.

Let $D$ be a countable and dense subset of $\mathrm{D}_{\mathcal{F}}^{\alpha}(P)$ and $x \in D$. There then exists $Q \in P^{\alpha}$ with $f(x) \leq \int f \,\mathrm{d}Q$ for all $f \in \mathcal{F}$.

Let us take the sequence $\{Q_n\}_{n \geq n_1}$ constructed in Proposition 9.

We split the set $\{n \in N : n \geq n_1\}$ into two sets, namely $N_Q$ (see expression (3)) and $\{n \in N : n \geq n_1\} \setminus N_Q$.

In Cascos and López-Díaz [3], it is proved that for the subsequence determined by $N_Q$, it holds that

$$\lim_{k} \left( \sup_{f \in \mathcal{F}} \left| \int f \,\mathrm{d}Q_{n_k} - \int f \,\mathrm{d}Q \right| \right) = 0 \qquad \text{a.s.}$$

Now, consider the subsequence determined by $\{n \in N : n \geq n_1\} \setminus N_Q$.

Given $n$ in this set, then $n_k < n < n_{k+1}$ for some $k$ with $n_k, n_{k+1} \in N_Q$. From (1), we have

$$Q_n(\cdot) = \frac{n_k}{n} Q_{n_k}(\cdot) + \frac{1}{n}(I_{(\cdot)}(X_{n_k+1}) + I_{(\cdot)}(X_{n_k+2}) + \cdots + I_{(\cdot)}(X_n)),$$

which implies that

$$\int f \,\mathrm{d}Q_n = \int \frac{n_k}{n} f \,\mathrm{d}Q_{n_k} + \frac{1}{n} \sum_{j=n_k+1}^{n} f(X_j).$$

Since $n_k \in N_Q$, equation (2) implies that

$$\int f \,\mathrm{d}Q_{n_k} = \int \frac{fg}{\|g\|_{L^1(P_{n_k})}} \,\mathrm{d}P_{n_k}.$$

Hence, we obtain that

$$\sup_{f\in\mathcal{F}}\left|\int f\,dQ_n - \int f\,dQ\right|$$

$$= \sup_{f\in\mathcal{F}}\left|\frac{n_k}{n}\int \frac{fg}{\|g\|_{L^1(P_{n_k})}}\,dP_{n_k} + \frac{1}{n}\sum_{j=n_k+1}^{n} f(X_j) - \int fg\,dP\right|$$

$$\leq \sup_{f\in\mathcal{F}}\left(\left|\frac{n_k}{n}\int \frac{fg}{\|g\|_{L^1(P_{n_k})}}\,dP_{n_k} - \int fg\,dP\right|\right) + \frac{1}{n}\sum_{j=n_k+1}^{n} f_{\mathcal{F}}(X_j).$$

When $n$ tends to $\infty$, since $n_k < n < n_{k+1}$, $n_k$ also tends to $\infty$ and, in accordance with (4), we obtain that

$$\lim_n \left(\sup_{f\in\mathcal{F}}\left|\frac{n_k}{n}\int \frac{fg}{\|g\|_{L^1(P_{n_k})}}\,dP_{n_k} - \int fg\,dP\right|\right)$$

$$= \lim_n \left(\sup_{f\in\mathcal{F}}\left|\frac{n_k}{n}\int f\,dQ_{n_k} - \int f\,dQ\right|\right) = 0 \quad \text{a.s.}$$

Note that $n_k/n_{k+1} \leq n_k/n \leq 1$ and from Proposition 8, it follows that $\lim_k n_k/n_{k+1} = 1$ a.s., so $\lim_k n_k/n = 1$ a.s.

Regarding the term $n^{-1}\sum_{j=n_k+1}^{n} f_{\mathcal{F}}(X_j)$, we trivially obtain that

$$\frac{1}{n}\sum_{j=n_k+1}^{n} f_{\mathcal{F}}(X_j) = \frac{1}{n}\sum_{j=1}^{n} f_{\mathcal{F}}(X_j) - \frac{1}{n}\sum_{j=1}^{n_k} f_{\mathcal{F}}(X_j)$$

$$= \frac{1}{n}\sum_{j=1}^{n} f_{\mathcal{F}}(X_j) - \frac{n_k}{n}\frac{1}{n_k}\sum_{j=1}^{n_k} f_{\mathcal{F}}(X_j).$$

On the other hand, $f_{\mathcal{F}} \in L^1(P)$ and since $\int f_{\mathcal{F}}\,dP = E(f_{\mathcal{F}}(X_1))$, the strong law of large numbers implies that

$$\lim_k \left(\frac{1}{n}\sum_{j=1}^{n} f_{\mathcal{F}}(X_j) - \frac{n_k}{n}\frac{1}{n_k}\sum_{j=1}^{n_k} f_{\mathcal{F}}(X_j)\right) = E(f_{\mathcal{F}}(X_1)) - E(f_{\mathcal{F}}(X_1)) = 0 \quad \text{a.s.}$$

As a consequence, we conclude that the sequence $\{Q_n\}_{n\geq n_1}$ satisfies

$$\lim_n \left(\sup_{f\in\mathcal{F}}\left|\int f\,dQ_n - \int f\,dQ\right|\right) = 0 \quad \text{a.s.}$$

So, we have a.s. that for all $\varepsilon > 0$, there exists $n_0 \in \mathbb{N}$ such that for all $n \geq n_0$, it holds that

$$\sup_{f\in\mathcal{F}}\left|\int f\,dQ_n - \int f\,dQ\right| \leq \frac{\varepsilon}{2}.$$

This implies that for all $f \in \mathcal{F}$, we have $f(x) - \varepsilon \leq \int f \, dQ_n$, so it holds a.s. that for all $\varepsilon > 0$, the point $x$ satisfies

$$x \in \liminf_n \mathrm{D}_{\mathcal{F}}^{\alpha,\varepsilon}(P_n).$$

Since the lower limit is closed and $D$ is countable, its topological closure satisfies

$$\mathrm{cl} D \subset \bigcap_{\varepsilon > 0} \liminf_n \mathrm{D}_{\mathcal{F}}^{\alpha,\varepsilon}(P_n) \qquad \text{a.s.}$$

and, finally,

$$\mathrm{D}_{\mathcal{F}}^{\alpha}(P) \subset \bigcap_{\varepsilon > 0} \liminf_n \mathrm{D}_{\mathcal{F}}^{\alpha,\varepsilon}(P_n) \qquad \text{a.s.},$$

which implies the desired result. $\square$

## Acknowledgements

The authors wish to thank Teófilo Brezmes from the Universidad de Oviedo for his suggestions and comments.

The financial support of the grant MTM2005-02254 of the Spanish Ministry of Education and Science is acknowledged.